\input amstex
\documentstyle{amsppt}

\magnification=\magstep1

\hoffset1 true pc \voffset2 true pc \hsize36 true pc \vsize50 true
pc

\footline{}

\catcode`\@=11 
\newfam\ssffam
\font\tenssf=cmss10 \font\eightssf=cmss8
\def\ssf{\fam\ssffam}
  \textfont\ssffam=\tenssf \scriptfont\ssffam=\eightssf
\def\sf#1{\leavevmode\skip@\lastskip\unskip\/%
       \ifdim\skip@=\z@\else\hskip\skip@\fi{\ssf#1}}
\def\rom#1{\leavevmode\skip@\lastskip\unskip\/%
        \ifdim\skip@=\z@\else\hskip\skip@\fi{\rm#1}}
\def\newmcodes@{\mathcode`\'"27\mathcode`\*"2A\mathcode`\."613A%
 \mathcode`\-"2D\mathcode`\/"2F\mathcode`\:"603A }
\def\operatorname#1{\mathop{\newmcodes@\kern\z@\fam\z@#1}}
\catcode`\@=\active

\catcode`\@=11
\def\raggedcenter@{\leftskip\z@ plus.3\hsize \rightskip\leftskip
 \parfillskip\z@ \parindent\z@ \spaceskip.3333em \xspaceskip.5em
 \pretolerance9999\tolerance9999 \exhyphenpenalty\@M
 \hyphenpenalty\@M \let\\\linebreak}
\catcode`\@=\active

\TagsOnRight
\define\m1{^{-1}}
\define\ov1{\overline}
\def\gp#1{\langle#1\rangle}
\def\ul2#1{\underline{\underline{#1}}}

\def\Y {\mathbin{\ssf Y@!@!@!@!}}
\def\char{\operatorname{char\, }}
\topmatter

\title
Normality in  group rings
\endtitle
\author
V.A.~Bovdi and S.~Siciliano
\endauthor
\dedicatory Dedicated to  Professor  P.M.~Gudivok on his  70th
birthday
\enddedicatory
\abstract Let $KG$ be the group ring of a group $G$ over  a
commutative ring $K$ with unity.  The rings $KG$ are described for
which $xx^\sigma=x^\sigma x$ for all $x=\sum_{g\in G}\alpha_gg\in
KG$, where \quad $x\mapsto x^\sigma=~\sum_{g\in
G}\alpha_gf(g)\sigma(g)$\quad  is an involution of $KG$; here $f:
G\to U(K)$ is  a homomorphism and $\sigma$ is an anti-automorphism
of order two of $G$.
\endabstract
\thanks
The research was supported by OTKA  No.T 037202 and  No.T 038059
\endthanks
\address
\hskip-\parindent {\rm  V.A.~Bovdi\newline Institute of
Mathematics, University of Debrecen\newline P.O.  Box 12, H-4010
Debrecen, Hungary\newline Institute of Mathematics and
Informatics, College of Ny\'\i regyh\'aza\newline S\'ost\'oi \'ut
31/b, H-4410 Ny\'\i regyh\'aza, Hungary
\newline
E-mail: vbovdi\@math.klte.hu}
\newline
\hskip-\parindent {\rm  S.~Siciliano
\newline
Dipartimento di Matematica ``E. De Giorgi", Universit\`{a} degli
Studi di Lecce\newline Via Provinciale Lecce-Arnesano,
73100-LECCE, Italy
\newline
E-mail: salvatore.siciliano\@unile.it}
\endaddress

\endtopmatter

\document
Let $R$ be a ring with unity. We denote by $U(R)$ the group of
units of $R$. A (bijective) map $\diamond: R \to R$ is called  an
{\it involution} if  for all $a, b\in R$ we have
$(a+b)^\diamond=a^\diamond+b^\diamond$, \quad
$(ab)^\diamond=b^\diamond\cdot a^\diamond$ \quad and \quad
$a^{\diamond^2}=a$. Let $KG$ be the group ring of a group $G$ over
the commutative ring $K$ with unity, let  $\sigma$ be an
antiautomorphism of order two  of $G$, and let  $f: G\to U(K)$ be a
homomorphism from  $G$ onto $U(K)$. For an element $x=\sum_{g\in
G}\alpha_gg \in KG$, we define $x^\sigma=~\sum_{g\in
G}\alpha_gf(g)\sigma(g) \in KG$. Clearly, $x\mapsto x^\sigma$ is
an involution of $KG$ if and only if \quad $g\sigma(g)\in Ker
f=\{h\in G\mid f(h)=1\}$\quad for all $g\in G$.

The ring $KG$ is said to be {\it $\sigma$-normal} if
$$
xx^\sigma=x^\sigma x \tag1
$$
for each $x\in KG$.
The properties  of the classical involution $x\mapsto x^{*}$
(where $*: g\mapsto g\m1$\quad  for  $g\in G$)  and the properties
of   normal group rings (i.e., $xx^*=x^*x$ for each $x\in KG$) have
been used  actively  for investigation of the group of units
$U(KG)$ of the group ring $KG$ (see \cite{1, 2}). Moreover, they also
have  important  applications in topology (see \cite{7, 8}).
Our aim is to describe the structure of the $\sigma$-normal group
ring  $KG$ for an arbitrary order $2$ antiautomorphism $\sigma$
of the group $G$. Note that  the descriptions of the classical
normal group rings and the twisted group rings were obtained in
\cite{1, 3} and \cite{4, 5} respectively.

The notation used throughout the paper is essentially standard.
\quad $C_n$  denotes the cyclic group of order $n$;\quad
$\zeta(G)$ and $C_G(H)$ are the center of the group $G$ and  the
centralizer  of $H$ in $G$,  respectively;\quad
$(g,h)=g^{-1}h^{-1}gh=g^{-1}g^h$\; ($g,h\in G$); \quad
$\gamma_{i}(G)$\quad is the $i${th} term of the lower central
series of $G$, i.e.,\quad  $\gamma_{1}(G)=G$ \quad and \quad
$\gamma_{i+1}(G)=\big(\gamma_{i}(G),G\big)$\quad for $i\geq 1$;
$\Phi(G)$ denotes the Frattini subgroup of $G$. We say that $G=A\Y
B$ is a central product of its subgroups $A$ and $B$ if $A$ and
$B$ commute elementwise and, taken together, they generate $G$, provided
that $A\cap B$ is a subgroup of $\zeta(G)$.

\smallskip
\noindent
 A non-Abelian $2$-generated nilpotent group  $G=\gp{a,b}$ with an
anti-automorphism $\sigma$ of order $2$ is  called a
$\sigma$-group if  $G'$ has order $2$, \; $\sigma(a)=a(a,b)$\quad
and \; $\sigma(b)=b(a,b)$.
\medskip

\noindent Our main result reads as follows.

\proclaim {Theorem} Let $KG$ be the noncommutative  group ring of
a group $G$ over a commutative ring $K$ and  $f: G\to U(K)$ \quad
 a homomorphism. Assume that  $\sigma$ is an antiautomorphism of
order two of $G$ such that $x\mapsto x^\sigma$ is an involution
of $KG$. Put $\frak{R}(G)=\{g\in G\mid \sigma(g)=g\}$. The group
ring   $KG$ is $\sigma$-normal\quad   if and only if\quad  $f:
G\to \{\pm1\}$,\quad $G$, $K$ and $\sigma$ satisfy one of  the
following conditions: \item{(i)}$G$ has an Abelian subgroup  $H$
of index $2$ such that \; $G=\gp{H,b}$,\quad $f(b)=-1$,\quad
$f(h)=1$,\quad $\sigma(b)=b$, \quad  and\quad
$\sigma(h)=b\m1hb=bhb\m1$ \quad for all\quad $h\in H$;
\item{(ii)}$G= H\Y\frak{C}$ is a central product of a
$\sigma$-group $H=\gp{a,b}$ and an abelian group $\frak{C}$ such
that   $G'= \gp{c\mid c^2=1}$ and $H\subset Ker(f)$. Moreover,
either  $\sigma(d)=d$ for all $d\in \frak{C}$,\quad
$\frak{C}\subset Ker(f)$,\quad  $\frak{R}(G)=\zeta(G)$, and
$$
G/\frak{R}(G)=\gp{a\zeta(G),\;b\zeta(G)}\cong C_2\times C_2
$$

or\quad   $\frak{R}(G)$ is of index $2$ in $\zeta(G)$ and
$$
G/\frak{R}(G)=\gp{g\frak{R}(G),\;h\frak{R}(G),\;d\frak{R}(G)}\cong C_2\times C_2\times C_2,
$$
where  $d\in \frak{C}$,\quad $\sigma(d)=dc$ and $f(d)=-1$;

\item{(iii)} $char(K)=2$,\quad   $G=S\Y\frak{C}$ is a central product of $S=\Y_{i=1}^{n} H_i$
and an Abelian group $\frak{C}$ such that $H_i=\gp{a_i,b_i}$ is a $\sigma$-group and   $G=Ker(f)$.
Moreover, $G'=\gp{c\mid c^2=1}$,\quad   $n\geq 2$, where $n$ is not necessarily a finite number,\quad
$\sigma(a_i)=a_ic$,\quad $\sigma(b_i)=b_ic$ \quad for all $i=1,2,\ldots$\quad  and\quad
$exp(G/\frak{R}(G))=2$.
\newline
Furthermore, if $n$ is finite, then either  $\sigma(d)=d$ for all $d\in \frak{C}$  and
$$
G/\frak{R}(G)=\times_{i=1}^{n} \gp{a_i\zeta(G),\; b_i\zeta(G)}\cong
\times_{i=1}^{2n} C_2,
$$
or\quad  $\frak{R}(G)$ is of index $2$ in $\zeta(G)$ and
$$
G/\frak{R}(G)=\times_{i=1}^{n} \gp{a_i\frak{R}(G),\; b_i\frak{R}(G)}\times \gp{d\frak{R}(G)}\cong
\times_{i=1}^{2n+1} C_2,
$$
where  $d\in \frak{C}$ and $\sigma(d)=dc$.
\endproclaim

\bigskip
Note that, in parts (i) and (ii) of the theorem, the group $\frak{C}$ may be equal to 1.

To make the statements less cumbersome, in what follows we will often speak of $\sigma$-normal group rings $KG$ without specifying the  homomorphism $f: G\to U(K)$ and
the anti-automorphism $\sigma$
 of order two of $G$. In order to prove the main theorem, we need some preliminary lemmas.

\proclaim {Lemma 1} Let $U(R)$ be the group of units of the ring
$R$, and let $x\mapsto x^\diamond$ be an involution of $R$. Suppose that
$xx^\diamond=x^\diamond x$ for all $x\in R$. If $a\in U(R)$, then
$a^\diamond=ta$,\quad  $at=ta$, and\quad  $t^\diamond=t\m1$, where
$t\in U(R)$.
\endproclaim

\demo{ Proof} Clearly $a^\diamond=at$ for some  $t\in U(R)$, and
$a^\diamond a=aa^\diamond$ implies  that  $ata=a^2t$ and $at=ta$.
Now
$a=a^{\diamond^2}=(at)^\diamond=t^{\diamond}at=t^{\diamond}ta$,\quad
whence $t^\diamond=t\m1$.{\text{\qed}}\enddemo

\proclaim {Lemma 2} Let $K$ be a commutative ring, let
$H=\gp{a,b}$ be a non-Abelian $2$-generated subgroup of a group $G$,
and let $f: G\to U(K)$ be a homomorphism. If the group ring  $KG$ is
$\sigma$-normal, then $f: H\to \{\pm 1\}$ and one of the following
conditions is fulfilled: \item{(i)}$f(a)=1$,\quad $f(b)=-1$, \quad
$\sigma(a)=(a,b)a$,\quad $\sigma(b)=b$,\quad $(b^2,a)=1$,\quad
$(ab)^2=(ba)^2$,\quad $((a,b),a)=1$, and \quad $((a,b),b)=(a,b)^{-2}$;
\item{(ii)}$f(a)=-1$,\quad $f(b)=1$, \quad $\sigma(a)=a$,\quad
$\sigma(b)=(a,b)b$,\quad $(b,a^2)=1$,\quad $(ab)^2=(ba)^2$,\quad
$((a,b),b)=1$, and \quad $((a,b),a)=(a,b)^{-2}$;
\item{(iii)}$f(a)=f(b)=-1$, \quad $\sigma(a)=a$,\quad
$\sigma(b)=b$,\quad $(a^2,b)=(a,b^2)=1$,\quad
$(ab)^2=(ba)^2$,and \quad $((a,b),ab)=1$; \item{(iv)} $f(a)=f(b)=1$,
\quad $\sigma(a)=(a,b)a$,\quad $\sigma(b)=(a,b)b$, \quad  and
\quad $\gp{a,b}$ is nilpotent of class $2$ and such that
$\gamma_2(\gp{a,b})$ of order $2$.
\endproclaim

\demo{ Proof} Let $KG$ be a $\sigma$-normal  ring. For any
noncommutative $a, b\in G$ we can put  $\sigma(a)=at$ and
$\sigma(b)=bs$, where $s,t\in G$. By Lemma 1,  $at=ta$,
$bs=sb$, $\sigma(t)=t\m1$, and $\sigma(s)=s\m1$. Set $x=a+b\in KG$.
Clearly, $x^\sigma=f(a)\sigma(a)+f(b)\sigma(b)$, and by (1) we have
$$
f(b)a\sigma(b)+ f(a)b\sigma(a) =f(a)\sigma(a)b+
f(b)\sigma(b)a.\tag2
$$
If\quad $a\sigma(b)=b\sigma(a)=\sigma(a)b=\sigma(b)a$,\quad then
we get $s=t$ and $ab=ba$,   a contradiction. Observe that if
three  of the  elements $\{a\sigma(b), b\sigma(a),\sigma(a)b,\sigma(b)a\}$
coincide, then $s=t$ and $ab=ba$, a contradiction. We consider
the following  cases:
\newline
{\bf 1}. \; $a\sigma(b)=b\sigma(a)$. By (2) it follows that
$$
f(a)+f(b)=0,\qquad a\sigma(b)=b\sigma(a),\qquad
\sigma(a)b=\sigma(b)a.\tag3
$$
{\bf 2}. \; $a\sigma(b)=\sigma(a)b$. This  yields $asb=atb$, so that
$s=t\in\zeta(H)$, and  (2) ensures
$$
f(a)=f(b),\qquad \sigma(a)=at,\qquad \sigma(b)=bt,\qquad
t\in\zeta(H). \tag4
$$
{\bf 3}. \;  $a\sigma(b)=\sigma(b)a$. Since
$b\sigma(b)=\sigma(b)b$ we  get   $\sigma(b)\in \zeta(H)$, a
contradiction.

Now put $x=a(1+b)$. Then\quad
$x^\sigma=(1+f(b)\sigma(b))f(a)\sigma(a)$\quad  and, by  (1),
$$
f(ab)a\sigma(ab)+f(a)ab\sigma(a)=f(a)\sigma(a)ab+f(ab)\sigma(ab)a.\tag5
$$
We shall treat the following  cases separately:
\newline
{\bf 1}. \;  $a\sigma(ab)=ab\sigma(a)$. Formula (5) implies that
$$
f(b)=-1,\qquad \sigma(b)=b,\qquad (\sigma(a)a)\cdot b=b\cdot
(\sigma(a)a). \tag6
$$
{\bf 2}. \;  $a\sigma(ab)=\sigma(a)ab$ and
$ab\sigma(a)=\sigma(ab)a$. By (1)  we have $ab=\sigma(b)a$
and\; $\sigma(b)=aba\m1$. Since $a\sigma(ab)=\sigma(a)ab$, we get
$aba\m1at=atb$ and $(b,t)=1$. Recall  that $\sigma(b)=bs=sb$. So,
by (5),
$$
f(b)=1,\qquad \sigma(b)=aba\m1,\qquad t\in\zeta(H), \quad
s=(a\m1,b\m1)=(b,a\m1). \tag7
$$
{\bf 3}. \; $a\sigma(b)\sigma(a)=\sigma(b)\sigma(a)a$. Then
$a\sigma(b)=\sigma(b)a$ and  $\sigma(b)\in \zeta(H)$, a
contradiction.

Assume that (3) and (6) are true. Then\quad  $f(b)=-1$,\quad $f(a)=1$,
\quad $\sigma(b)=b$,\quad $\sigma(a)=b\m1ab=bab\m1$, whence
$(b^2,a)=1$. Since\quad
$\sigma(a)=a(a\m1b\m1ab)=(bab\m1a\m1)a$,\quad  we get \quad
$a\m1b\m1ab=bab\m1a\m1$\quad  and\quad  $(ab)^2=(ba)^2$.
Obviously,
$$
b\m1(a,b)b=b\m1(bab\m1a\m1)b=ab\m1a\m1b=aba\m1b\m1=(a,b)\m1,
$$
so that \quad $((a,b),b)=(a,b)^{-2}$, \quad and statement (i) of our Lemma follows.

If  (3) and (7) are fulfilled, then $f(b)=1$,  $f(a)=-1$,
$\sigma(a)=a$,   $\sigma(b)=a\m1ba$,  and
$(a^2,b)=1$.  Since\quad
$\sigma(b)=b(b\m1a\m1ba)=(aba\m1b\m1)b$,\quad  we obtain \quad
$s=b\m1a\m1ba=aba\m1b\m1$\quad  and\quad  $(ab)^2=(ba)^2$.
Therefore,  $\sigma(a)=a$,\quad  $\sigma(b)=a\m1ba$,\quad
$(a^2,b)=1$, and  we arrive at statement (ii).

Assume (4) and (6). Then $f(a)=f(b)=-1$, \quad
$\sigma(a)=a$,\quad   $\sigma(b)=b$, and $(a^2,b)=1$. Moreover
$f(ab)=1$ and \quad $\sigma(ab)=ba=a\m1(ab)a$. We put $x=b(1+a)$.
Clearly,\quad  $x^\sigma=(a-1)b,$\;  and (1) implies\quad
$b^2a+b^2a^2=ab^2a+ab^2$,\quad whence \; $(b^2,a)=1$.
Thus, statement (iii) of our lemma is fulfilled.

Finally, if (4) and (7) are true, then $f(a)=f(b)=1$ and
$(b,a\m1)\in \zeta(H)$. Using the identity
$(\alpha\beta,\gamma)=(\alpha,\gamma) (\alpha,\gamma,\beta)
(\beta,\gamma)$, where $\alpha,\beta,\gamma\in G$,  we see
that $1=(a\m1a,b)=(a\m1,b)(a,b)$,   whence $s= (b,a\m1)=(a,b)\in
\zeta(H)$ and $\sigma(a)=(a,b)a$,\quad $\sigma(b)=(a,b)b$. Since
$a=\sigma^2(a)$ and $(a,b)\in \zeta(H)$, we have that
$(a,b)^2=1$, which yields statement (iv). The proof is
complete. {\text{\qed}}
\enddemo

\proclaim {Lemma 3} Let $KG$ be a $\sigma$-normal group ring of a
non-Abelian group $G$. Then $H=\gp{w\in G\mid \sigma(w)\not=w}$ is
a normal subgroup in $G$. If $H$ is Abelian,  then $G$ satisfies
 statement (i) of the Theorem.
\endproclaim

\demo{Proof} Set $W=\{w\in G\mid \sigma(w)\not=w\}$. Let $g\not\in
W$ be such that  $g^2\not\in\zeta(G)$. Then $(g^2,h)\not=1$ and
$(g,h)\not=1$, respectively,   for some $h\in G$.

We consider the following cases. \item{{\bf 1}.} $\char(K)\not=2$.
Since $\sigma(g^2)=g^2$, we can use Lemma 2 for  the group $\gp{g^2,h}$
to show that \quad  $-1=f(g^2)=(\pm 1)^2=1$, a contradiction.
\item{{\bf 2}.} $\char(K)=2$. Using Lemma 2 for  $\gp{g,h}$, we
get $(g^2,h)=1$, again  a contradiction.

Thus,  $g^2\in\zeta(G)$ for any $g\not\in W$. Now, if  $w\in W$,
$g\in G\setminus W$ and $g\m1wg\not\in W$, then
$\sigma(g\m1wg)=g\m1wg$\quad  and
$$
g\m1wg=\sigma(g\m1wg)=g\sigma(w)g^{-2}g=g\m1\sigma(w)g
$$
so that $\sigma(w)=w$, a contradiction. Therefore, $g\m1wg\in W$ and
the subgroup $H=\gp{W}$ is normal in $G$.

Suppose that $H=\gp{W}$ is Abelian. If $a\in W$ and $c\in C_G(W)\setminus H$, then
$ca\not\in H$. Therefore, $ca=\sigma(ca)=\sigma(a)c$,  whence
\quad $\sigma(a)=a$, a contradiction. This shows that  $C_G(W)=H$ and
for each   $b\not\in H$ there exists $w\in W$ such that   $(b,w)\not=1$.

We claim that   if $b_1,b_2\in G\setminus H$, then $b_1b_2\in H$.  The following cases will be treated separately:
\newline
{\bf 1.} \; $char(K)\not=2$ and    $b_1b_2\in G\setminus H$. For each
$b_i$ we choose  $w_i\in W$ such that $(b_i,w_i)\not=1$. By (i) or
(ii) of Lemma 2, in $\gp{w_i,b_i}$ we have  $f(b_i)=-1$, so that
$f(b_1b_2)=1$ and there exists  $w\in W$, for which
$(b_1b_2,w)\not=1$. Since $\sigma(b_1b_2)=b_1b_2$,\quad   by (i)
or (ii) of  Lemma 2 we get $f(b_1b_2)=-1$, a contradiction.
\newline
{\bf 2.} \; $char(K)=2$ and   $b_1b_2\in G\setminus H$. Obviously,
$b_1b_2=\sigma(b_1b_2)=b_2b_1$, whence $(b_1,b_2)=1$. Now, there is
$w\in W$ with $(w,b_1)\not=1$, and by Lemma 2 we get
$\sigma(w)=b_1\m1wb_1=b_1wb_1\m1$. Furthermore, $b_1b_2w\in
G\setminus H$ and
$$
b_1b_2w=\sigma(b_1b_2w)=\sigma(w)b_1b_2=b_1wb_2,
$$
implying $(b_2,w)=1$. Now $(b_1,b_2w)=(b_1,w)\not=1$ and $b_2w\in
G\setminus H$; applying Lemma 2 in $\gp{b_1,b_2w}$, we obtain
$b_2w=\sigma(b_2w)=\sigma(w)b_2$ and $\sigma(w)=w$, a
contradiction.

We have proved  that $b_1b_2\in H$ for every  $b_1,b_2\in G\setminus
H$. Therefore, $G=\gp{H,b\mid b\not\in H, b^2\in H}$,\quad $f(b)=-1$,
and $f(h)=1$ \quad  for all $h\in H$.

Finally, let $w\in W$ be such that $(b,w)=1$. Since $b\not\in H=C_G(W)$, there exists $w_1\in W$  with $w_1\not=b\m1w_1b=\sigma(w_1)$.  Clearly, we have $(ww_1,b)\not=1$;
using  Lemma 2 for $\gp{ww_1,b}$, we obtain
$$
\sigma(w)\sigma(w_1)=\sigma(w_1w)=b\m1w_1wb=\sigma(w_1)w,
$$
whene $\sigma(w)=w$, a contradiction. Thus, \quad
$b\m1hb=\sigma(h)$ \quad  for all $h\in H$.{\text{\qed}}
\enddemo

\proclaim {Lemma 4} Let $KG$ be a $\sigma$-normal group ring,
let $W=\{w\in G\mid \sigma(w)\not=w\}$, and let $a,b\in W$ be such that
$(a,b)\not=1$. Put \quad $\frak{R}=\{g\in G\mid \sigma(g)=g\}$ and
$\frak{C}=C_G(\gp{a,b})$. Then  $\gp{a,b}$ is a
$\sigma$-group,\qquad $\Phi(\gp{a,b})=\zeta(\gp{a,b})=\{g\in
\gp{a,b}\mid \sigma(g)=g\}$,\qquad and
$$
\sigma(g)=\cases g& \text{if}\qquad  g\in\zeta(\gp{a,b});\\
g(a,b)& \text{if}\qquad  g\not\in \zeta(\gp{a,b}).
\endcases
$$
Moreover,  $G=\gp{a,b}\Y \frak{C}$, and either  $\sigma(c)=(a,b)c$
or   $\sigma(c)=c$ where $c\in \frak{C}$. Also,  the
following is true: \item{(i)} if $\frak{C}$ is abelian, then  $G$
satisfies statement (ii) of the Theorem; \item{(ii)} if $
\frak{C}$ is not Abelian, then $char(K)=2$.
\endproclaim

\demo{Proof} Let $a,b\in W$ satisfy $(a,b)\not=1$. By  Lemma 2,
$f(a)=f(b)=1$, $\gp{a,b}$ is nilpotent of class $2$ and such that $|\gamma_2(\gp{a,b})|=2$, and \quad
$\sigma(a)=b\m1ab$ and $\sigma(b)=a\m1ba$. Thus $\gp{a,b}$ is a $\sigma$-group. Any element
$g\in \gp{a,b}$ can be written  as $g=a^ib^j(a,b)^k$,
where $i,j,k\in \Bbb N$. Since $\sigma(g)=g$, we conclude  that $i$ and $j$
are even. Now by [6, Theorem 10.4.1 and  Theorem 10.4.3] we
obtain \quad
$$
\Phi(\gp{a,b})=\zeta(\gp{a,b})=\{g\in
\gp{a,b}\mid \sigma(g)=g\}.
$$
Suppose $c\in W$ and $(a,c)\not=1$. Again by Lemma 2, $\gp{a,c}$ is
nilpotent of class $2$ and\quad $\sigma(a)=c\m1ac=b\m1ab$, so that
\quad $(a,b)=(a,c)$. Now, let\quad $c,d\in W$ be such that
$(c,d)\not=1$ and $\gp{c,d}\in \frak{C}$.\quad Obviously,\quad
$(ac,b)=(a,b)\not=1$\quad and\quad $(ac,d)=(c,d)\not=1$.\quad  By
Lemma 2,  $\sigma(ac)=b\m1(ac)b=d\m1(ac)d$\quad and\quad
$(a,b)=(c,d)$, which shows that \quad $H'$ has order two  and is central in
$G$.

Let  $g\in G\setminus \frak{C}\cdot \gp{a,b}$. If $(a,g)\not=1$,
then using  Lemma 2 for $\gp{a,g}$, we get
$\sigma(a)=g\m1ag=b\m1ab$ and $(a,g)=(a,b)$. Similarly, if
$(b,g)\not=1$, then $(b,g)=(a,b)$.

The following cases are possible:
\item{{\bf 1.}}  $(g,a)=1$ and  $(g,b)\not=1$. Then we have
$(ga,b)=(ga,a)=1$, whence  $g=(ga)\cdot
a\m1\in \frak{C}\cdot \gp{a,b}$.
\item{{\bf 2.}} $(g,a)\not=1$ and  $(g,b)=1$. Then we have
$(gb,a)=(gb,b)=1$, whence $g=(gb)\cdot b\m1\in
\frak{C}\cdot \gp{a,b}$.
\item{{\bf 3.}} $(g,a)\not=1$ and $(g,b)\not=1$. Then we have
$(gab,b)=(gab,a)=1$, whence $g=(gab)\cdot
(ab)\m1\in \frak{C}\cdot \gp{a,b}$.

Since each of these  cases leads to  a contradiction, we have $G=\frak{C}\Y \gp{a,b}$.

Let $d\in \frak{C}\setminus H$. Since  $\sigma(ad)=ad$ we get
$ad=\sigma(ad)=\sigma(a)d$, whence $\sigma(a)=a$, a contradiction.
Since $G=\frak{C}\cdot \gp{a,b}$, it follows  that $G=H=\gp{W}$.
If $d\in \zeta(G)\cap W$, then  $\sigma(ad)=ad$, and\quad
$(ad,b)=(a,b)\not=1$; using Lemma 2 for $\gp{ad,b}$, we obtain
$$
-1=f(ad)=f(a)f(d)=f(d).
$$
Now, we let $\zeta(G)\cap W=\emptyset$ and put $x=ac+b$, where $c\in
\frak{C}$. Then there exists $d\in G$ such that $(c,d)\not=1$ and
 Lemma 2 implies that $f(g)=1$ for all $g\in G$.\quad Thus,
$x^\sigma=(a\sigma(c)+b)(a,b)$,\quad and by (1) we have\quad
$(\sigma(c)-c)(1-(a,b))=0$. It follows that either
$\sigma(c)=c$,\quad  or\quad $char(K)=2$ and $\sigma(c)=(a,b)c$.
Therefore, if $\frak{C}$ is Abelian, we obtain statement
(ii) of the Theorem.

Finally,  assume that $char(K)\not=2$. Suppose that
there exist $c,d\in \frak{C}$ such that  $(c,d)\not=1$. If
$\sigma(c)=c$, then $f(c)=1$ by what has already been proved,
but, by Lemma 2 in $\gp{c,d}$, we have $f(c)=-1$, a contradiction.
Therefore,  $c \in W$ and similarly $d \in W$. We put $x=ac+d$.
Clearly, $x^\sigma=ac+d(a,b)$ and $(a,b)=1$ by (1), a
contradiction. Thus, if $\frak{C}$ is not Abelian, then
$char(K)=2$, and the proof is complete. {\text{\qed}}
\enddemo
\smallskip
Now we are in a position to prove our main theorem.
 \smallskip

\noindent {\it Proof of the ``if" part of the Theorem.} Set
$W=\{w\in G\mid \sigma(w)\not=w\}$ and $H=\gp{W}$. If $H$ is
Abelian, then,  statement (i) of the
Theorem is valid for $G$.

Suppose that $H$ is non-Abelian and that $a,b\in W$  satisfy
$(a,b)\not=1$. By Lemma 4, $G=\gp{a,b}\Y \frak{C}=\gp{W}$, where
$\frak{C}=C_G(\gp{a,b})$. If $\frak{C}$ is Abelian, then statement (ii) of our Theorem
is valid for $G$ by Lemma 4.

Let  $c,d\in C_G(\gp{a,b})$ such that  $(c,d)\not=1$ (i.e.
$\frak{C}$ is non-Abelian). By Lemma 4, we have  $char(K)=2$. If
$c,d\in W$, then, by Lemma 4,
$$
G=C_G(\gp{a,b})\cdot \gp{a,b}=C_G(\gp{c,d})\cdot \gp{c,d}.
$$
Obviously, $C_G(\gp{a,b})\cap \gp{a,b}\subseteq \zeta(G)$. Therefore, $G$
contains the subgroup  $H_2=\gp{a,b}\Y \gp{c,d}$, which  cannot
be a direct product because $G'$ has  order $2$.

Since $G'\subseteq \frak{R}(G)$, we see that   $G/\frak{R}(G)$ is
an elementary Abelian $2$-group. Let $\tau: G\to
G/\frak{R}(G)=\times_{i\geq 1}\gp{a_i\mid a_i^2=1}$, such that
$\tau^{-1}(a_1)=a$ and $\tau^{-1}(a_2)=b$. We put
$\overline{a}_i=\tau^{-1}(a_i)$ for all $i\geq 3$ and
$\frak{B}=\{a_i\mid i\geq 3\}$.

Suppose that for some $s\geq 3$ we have
$(\overline{a}_s,\overline{a}_i)=1$ for all $i\geq 3$. Such an
element is unique, because if $\overline{a}_t\not=\overline{a}_s$
commutes  with all $\overline{a}_s$, then
$\sigma(\overline{a}_t\overline{a}_s)=\overline{a}_s\overline{a}_t$,
whence $a_sa_t=1$, a contradiction. Put $\frak{B}=\frak{B}\setminus
a_s$, $b_0=a_s$, $b_1=a_1$ and $b_2=a_2$. Note that if such an element
$a_s$ does not exist, then  we put $b_0=1$.

Choose $a_i\in \frak{B}$. There is $a_j\in \frak{B}$ such that
$(\overline{a}_i,\overline{a}_j)\not=1$, and we consider the
following cases.

{\bf 1.} $\overline{a}_i,\overline{a}_j\in W$. Put $b_3=a_i$,
$b_4=a_j$ and $\frak{B}=\frak{B}\setminus \{a_i,a_j\}$.

{\bf 2.}   $\overline{a}_i\in W$ and $\overline{a}_j\not\in W$.
Clearly,
$\gp{\overline{a}_1,\overline{a}_2}\Y\gp{\overline{a}_i,\overline{a}_j}\cong
\gp{\overline{a}_1\overline{a}_i,\overline{a}_2}\Y\gp{\overline{a}_i,\overline{a}_2\overline{a}_j}$
and
$\overline{a}_1\overline{a}_i,\overline{a}_2,\overline{a}_i,\overline{a}_2\overline{a}_j\in
W$. Put $b_1=\tau(\overline{a}_1\overline{a}_i)$, $b_2= a_2$,
$b_3=a_i$, $b_4=\tau(\overline{a}_2\overline{a}_j)$ and
$\frak{B}=\frak{B}\cup \{a_1,a_2\}\setminus \{b_1,b_2, b_3,b_4\}$.

{\bf 3.}  $\overline{a}_i, \overline{a}_j\not\in W$. Obviously, we have
$\overline{a}_i\overline{a}_j\in W$, so that this case reduces
to the preceding one.

Furthermore, if $C_G(\gp{\overline{b}_1,\overline{b}_2}\Y
\gp{\overline{b}_3,\overline{b}_4})$ contains a noncommuting pair
of elements, then  this pair can be chosen in $W$. By  continuing   this process we can conclude  that $G$
contains a subgroup ${\frak M}=A_1\Y A_2\Y\cdots$ that is a
central product, where each $A_i=\gp{\; g_i,h_i\;}$ is a
$\sigma$-group, and $C_G({\frak M})$ is Abelian. Applying
Lemma 4, we arrive at statement (iii) of the Theorem, and the proof
is done.

{\it Proof of the ``only if" part of the Theorem.}
(i) We can write  any $x\in KG$ as
$x=x_1+x_2b$, where $x_i\in KH$. Clearly, \quad
$x^\sigma=x_1^\sigma+f(b)\sigma(b)x_2^\sigma=x_1^\sigma-x_2b$\quad
and
$$
xx^\sigma=x_1x_1^\sigma-x_2x_2^{\sigma}b^2=x_1^{\sigma}x_1-x_2^{\sigma}x_2b^2=x^{\sigma}x,
$$
so that $KG$ is a $\sigma$-normal ring.

(ii) Any $x\in KH$ can be written  as \quad
$x=x_0+x_1g+x_2h+x_3gh$,\quad where $x_i\in K\gp{g^2, h^2, c}$ and
$c=(g, h)$. Clearly,\quad $x^\sigma=x_0+(x_1g+x_2h+x_3gh)c$\quad
and\quad $xx^\sigma=x^{\sigma}x$,\quad   so that $KH$ is
$\sigma$-normal. Suppose that $\sigma(d)=dc$, with $c=(a,b)$.  Any
$x\in KG$ can be written as $x=(w_0+u_1)+(w_2+u_3)d$,\quad
where\quad $u_1=\alpha_1a+\alpha_2b+\alpha_3ab$,\quad
$u_3=\beta_1a+\beta_2b+\beta_3ab$,\quad and\quad $\alpha_i,\beta_i,
w_0,w_2\in K\frak{R}$. Then $x^\sigma=(w_0+u_1c)-
(w_2+u_3c)dc$\quad and\quad
$xx^{\sigma}-x^{\sigma}x=(u_3u_1-u_1u_3)(1+c)d$. Since
$ab-ba=ba(c-1)$ and $c^2=1$, it follows that
$xx^{\sigma}-x^{\sigma}x=0$.   Thus, $KG$ is $\sigma$-normal. In the case where
 $\sigma(d)=d$ the proof is similar.

(iii) Put $G_n=A_1\Y\cdots\Y
A_{n}$, where $A_i=\gp{a_i,b_i\mid c=(a_i,b_i)}$ is a $\sigma$-subgroup. We  use
induction on $n$. Any $x\in KG_n$ can be written as\quad
$x=x_0+x_1a_n+x_2b_n+x_3a_nb_n$, \quad where $x_i\in K\gp{G_{n-1},
a_n^2, b_n^2}$. Obviously, \quad
$x^\sigma=x_0^\sigma+(x_1^{\sigma}a_n+x_2^\sigma
b_n+x_3^{\sigma}a_nb_n)c$.  Since $KG_{n-1}$ is
$\sigma$-normal, we get
$x_ix_i^{\sigma}=x_i^{\sigma}x_i$   and
$x_i^{\sigma}(1+c)=x_i(1+c)$. The formula
$$
(x_i+x_j)(x_i+x_j)^{\sigma}=(x_i+x_j)^{\sigma}(x_i+x_j)
$$
shows   that
$$
x_ix_j^{\sigma}+x_jx_i^{\sigma}=x_i^{\sigma}x_j+x_j^{\sigma}x_i.
$$
Proceeding as in the preceding case, we conclude that
$$
xx^{\sigma}=x^{\sigma}x,
$$ and the proof is complete.{\text{\qed}}


\bigskip
\bigskip

\Refs

\ref\no 1\by S.\, D. Berman \paper On the equation $x\sp m=1$ in
an integral group ring. (Russian) \jour Ukrain. Mat. \v Z. \vol
7(3) \yr 1955 \pages 353--261
\endref

\ref\no 2\by A.\, A. Bovdi \paper Unitarity of the multiplicative
group of an integral group ring. (Russian) \jour  Mat. Sb. (N.S.)
\vol 119(161) \yr 1982 \pages 387--400
\endref

\ref\no 3\by A.\, A. Bovdi, P.\, M. Gudivok, M.\, S. Semirot
\paper Normal group rings \jour Ukrain. Mat. Zh. \vol 37 \yr 1985
\pages 3--8
\endref


\ref\no 4 \by V.\, A. Bovdi \paper Normal twisted group rings.
(Russian) \jour Dokl. Akad. Nauk Ukrain. SSR Ser. A \yr 1990\vol
7\pages  6--8
\endref

\ref\no 5 \by V.\, A. Bovdi \paper Structure of normal twisted
group rings \jour Publ. Math. Debrecen \vol 51 \yr 1997 \pages
279--293
\endref





\ref\no 6 \by M. Hall \book Group theory \publ Springer \publaddr
Berlin \yr 1971, 410 p.
\endref


\ref\no 7 \by M.A. Knus, A. Merkurjev, M.Rost, J.-P. Tignol
\book The book of involutions
\publ AMS Colloquium Publications, AMS Providence, 44
\yr 1998 \pages 593
\endref


\ref\no 8 \by S.P. Novikov  \paper  Algebraic construction and
properties of Hermitian analogs of $K$-theory over rings with
involution from the viewpoint of Hamiltonian formalism.
Applications to differential topology and the theory of
characteristic classes.I. II \jour Math. USSR-Izv. \vol 4 \yr 1970
\pages 257--292
\endref


\endRefs
\enddocument